\documentclass[12pt]{article}

\usepackage{amssymb, amsmath, amsfonts, amsthm, mathrsfs, latexsym, epsfig,  psfrag, graphicx}
\usepackage[all]{xy}
\xyoption{v2} \xyoption{knot}\xyoption{arc} \xyoption{poly}
\xyoption{curve}

\usepackage{amsmath,amsthm,amsfonts,amssymb,amscd}
\def\qed{{\unskip\nobreak\hfil\penalty50
\hskip2em\hbox{}\nobreak\hfil$\square$
\parfillskip=0pt \finalhyphendemerits=0\par}\medskip}

\def\prf{\trivlist \item[\hskip \labelsep{\bf Proof\ }]}

\def\ro{\rho}
\def\ad{{v_0}}
\def\ra{\rangle}

\def\C{{\Bbb {C}}}
\def\End{{\mathrm {End}}}

\def\Hom{{\mathrm {Hom}}}

\def\col{{\mathrm{col}}}
\def\a{\alpha}

\def\e{\varepsilon}

\def\la{\lambda}

\def\om{\omega}
\def\Om{\Omega}

\def\Diff{{\mathrm {Diff}}}
\def\Mob{{\rm\textsf{M\"ob}}}

\def\End{{\mathrm {End}}}
\def\Hom{{\mathrm {Hom}}}

\def\a{\alpha}

\def\e{\varepsilon}
\def\r{\rho}

\def\la{\lambda}

\def\Om{\Omega}

\def\ro{{\rho}}

\newtheorem{theorem}{Theorem}[section]
\newtheorem{lemma}[theorem]{Lemma}
\newtheorem{conjecture}[theorem]{Conjecture}
\newtheorem{corollary}[theorem]{Corollary}
\newtheorem{definition}[theorem]{Definition}

\newtheorem{proposition}[theorem]{Proposition}
\newtheorem{remark}[theorem]{Remark}

\def\Diff{{\mathrm {Diff}}}
\def\Hom{{\mathrm{Hom}}}
\def\Mob{{\rm\textsf{M\"ob}}}

\def\res{\!\restriction\!}
\def\A{{\cal A}}

\def\B{{\cal B}}
\def\C{{\cal C}}

\def\I{{\cal I}}

\def\H{{\cal H}}

\def\Z{{\mathbb Z}}

\renewcommand{\qed}{\ \hfill $\blacksquare$}

\newcommand{\bdef}{\begin{definition}}
\newcommand{\blem}{\begin{lemma}}
\newcommand{\bprop}{\begin{proposition}}
\newcommand{\bthm}{\begin{theorem}}
\newcommand{\bcor}{\begin{corollary}}
\newcommand{\bconj}{\begin{conjecture}}
\newcommand{\ede}{\end{definition}}
\newcommand{\elem}{\end{lemma}}
\newcommand{\eprop}{\end{proposition}}
\newcommand{\ethm}{\end{theorem}}
\newcommand{\ecor}{\end{corollary}}
\newcommand{\econj}{\end{conjecture}}
\newcommand{\brem}{\begin{remark}}
\newcommand{\erem}{\end{remark}}

\newcommand{\ba}{\begin{array}}
\newcommand{\ea}{\end{array}}
\newcommand{\bea}{\begin{eqnarray}}

\renewcommand{\mod}{\mbox{mod}}

\title{\huge On examples of intermediate subfactors from conformal field theory \\}
\author{
{\sc Feng Xu}\footnote{Supported in part by NSF grant  and an academic senate grant from UCR.}\\
Department of Mathematics\\
University of California at Riverside\\
Riverside, CA 92521\\
E-mail: {\tt xufeng@math.ucr.edu}}

\begin{document}
\date{}
\maketitle

\begin{abstract}
Motivated by our subfactor generalization of Wall's conjecture,  in
this paper we determine all  intermediate subfactors for conformal
subnets corresponding to four infinite series of conformal
inclusions, and as a consequence we verify that these series of
subfactors verify our conjecture. Our results can be stated in the
framework of Vertex Operator Algebras. We also verify our conjecture
for Jones-Wassermann subfactors from representations of Loop groups
extending our earlier results.

\end{abstract}

\newpage

\section{Introduction}

Let $M$ be a factor represented on a Hilbert space and $N$ a
subfactor of $M$ which is irreducible, i.e.,$N'\cap M= \mathbb{C}$.
Let $K$ be an intermediate von Neumann subalgebra for the inclusion
$N\subset M.$ Note that $K'\cap K\subset N'\cap M = \mathbb{C},$ $K$
is automatically a factor. Hence the set of all intermediate
subfactors for $N\subset M$ forms a lattice under two natural
operations $\wedge$ and $\vee$ defined by:
\[
K_1\wedge K_2= K_1\cap K_2, K_1\vee K_2= (K_1\cup K_2)''.
\]
The commutant map $K\rightarrow K'$ maps an intermediate subfactor
$N\subset K\subset M$ to $M'\subset K'\subset N'.$ This map
exchanges the two natural operations defined above.

Let $M\subset M_1$ be the Jones basic construction of $N\subset M.$
Then  $M\subset M_1$ is canonically isomorphic to $M'\subset N'$,
and the lattice of intermediate subfactors for $N\subset M$  is
related to the lattice of intermediate subfactors for $M\subset M_1$
by the commutant map defined as above.

 Let $G_1$ be a group and $G_2$ be a subgroup of $G_1$. An
interval sublattice $[G_1/G_2]$ is the lattice formed by all
intermediate subgroups $K, G_2\subseteq K\subseteq G_1.$

By cross product construction and Galois correspondence,  every
interval sublattice of finite groups can be realized as intermediate
subfactor lattice of finite index. Hence the study of intermediate
subfactor lattice of finite index is a natural generalization of the
study of interval sublattice of finite groups. The study of
intermediate subfactors has been very active in recent years(cf.
\cite{BJ},\cite{GJ}, \cite{JXu}, \cite{ILP}, \cite{Longo4}, and
\cite{X-adv} for a partial list).

In 1961 G. E. Wall conjectured that the number of maximal subgroups
of a finite group $G$ is less than $|G|$, the order of $G$ (cf.
\cite{wall61}). In the same paper he proved his conjecture when $G$
is solvable. See \cite{lie} for more recent result on Wall's
conjecture.

Wall's conjecture can be naturally generalized to a conjecture about
maximal elements in the lattice of intermediate subfactors. What we
mean by maximal elements are those subfactors $K\neq M,N $ with the
property that if $K_1$ is an intermediate subfactor and $K\subset
K_1,$ then $K_1=M$ or $K.$ Minimal elements are defined similarly
where $N$ is not considered as an minimal element. When $M$ is the
cross product of $N$ by a finite group $G$, the maximal elements
correspond to maximal subgroups of $G,$ and the order of $G$ is the
dimension of second higher relative commutant. Hence a natural
generalization of Wall's conjecture as proposed in \cite{X-ghj} is
the following:
\begin{conjecture}\label{wall}
Let $N\subset M$ be an irreducible subfactor with finite index. Then
the number of maximal  intermediate subfactors is less than
dimension of $N'\cap M_1$ (the dimension of second higher relative
commutant of $N\subset M$).
\end{conjecture}
We note that since maximal intermediate subfactors in $N\subset M$
correspond to minimal intermediate subfactors in $M\subset M_1,$ and
the dimension of second higher relative commutant remains the same,
the conjecture  is equivalent to a similar conjecture as above with
maximal replaced by minimal.\par

In \cite{X-ghj},\cite{GX}, Conjecture \ref{wall} is verified for
subfactors coming from certain conformal field theories and
subfactors which are more closely related to groups and more
generally Hopf algebras. In this paper we investigate Conjecture
\ref{wall} for conformal subnets $\A\subset \B$ (cf. Definition
\ref{ext}) with finite index. Then Conjecture \ref{wall} in this
case states:
\begin{conjecture}\label{netwall}
Suppose that conformal subnets $\A\subset \B$ (cf. Definition
\ref{ext}) has finite index. Then the number of minimal (resp.
maximal) subnets between $\A$ and $\B$ is less than the dimension of
the space of bounded maps from the vacuum representation of $\B$ to
itself which commutes with the action of $\A$.
\end{conjecture}
In the above conjecture we have included both maximal and minimal
cases since the dual of conformal subnet $\A\subset \B$ is not
conformal subnet. It is also straightforward to phrase the above
conjecture in terms of Vertex Operator Algebras (VOAs) and its
sub-VOAs.\par

Note that any finite group $G$ is embedded in a finite symmetric
group $S_n,$ and using the theory of permutation orbifolds as in
\cite{X-permu} we can always find a completely rational net $\B$
such that $G$ acts properly on $\B$ and with fixed point subnet
$\A.$ In this case the intermediate subnets between $\A$ and $\B$
are in one to one correspondence with subgroups of $G.$ So in this
orbifold case the minimal version of  Conjecture \ref{netwall} is
equivalent to Wall's conjecture. Hence Conjecture \ref{netwall} is
highly nontrivial even if we assume that $\B$ is completely
rational. \par Though the orbifold case of Conjecture \ref{netwall}
in general is out of reach at present, there are very interesting
other examples of  subnets coming from conformal field theory (CFT).
A large class of such examples come from conformal inclusions (cf.
\S\ref{ci}), and they provide a large class of subfactors which are
not related to groups. In view of Conjecture \ref{netwall} it is a
natural question to investigate intermediate subnets of such
examples, and this is the main goal of our paper.\par

Our  results Th. \ref{c123}, Th. \ref{normalcase} give a complete
list of intermediate conformal subnets in subnets coming from four
infinite series of maximal conformal inclusions, and as consequence,
we are able to verify Conjecture \ref{netwall} in these examples.
Our results show that the intermediate subnets in these examples are
very rare. The key idea behind the proof of Th. \ref{c123} is the
property of induced adjoint representation:  Prop. \ref{key} shows
that such induced representation is always irreducible when the
intermediate subnet does not have additional weight $1$ element. By
locality consideration in Lemma  \ref{loc2}  this forces the
intermediate subnet to be simply simple current extensions when it
has no additional weight $1$ element. In the case when the
intermediate net has additional weight $1$ element, we use smeared
vertex operators as in \cite{X-coset} and maximality of conformal
inclusions to show that the intermediate subnet is in fact the
largest net. The proof makes use of  the analogue of statement in
VOA theory that weight $1$ element of a VOA forms a Lie algebra. The
proof of Th. \ref{normalcase} is much simpler and make use of normal
inclusions as in \S4.2  of \cite{X-m}.
\par By using properties of smeared vertex operators in
\S\ref{smear}, we can translate our results in Th. \ref{c123}, Th.
\ref{normalcase} into statements about intermediate VOAs (cf. Th.
\ref{voastatement}).  We think  it is an interesting question to
find a VOA proof of Th. \ref{voastatement}.\par

In \S4  we extend our earlier results in \S5 of \cite{X-adv} on
Jones-Wassermann subfactors and we verify that these subfactors
verify Conjecture \ref{wall}.\par

In addition to what are already described as above, we have included
a preliminary section \S2 where we introduce the basic notion of
conformal nets, subnets , conformal inclusions, and induction to
describe the background of our results in \S3 and \S4.

\section{Preliminaries}
\subsection{Preliminaries on sectors}

Given an infinite factor $M$, the {\it sectors of $M$}  are given by
$$\text{Sect}(M) = \text{End}(M)/\text{Inn}(M),$$
namely $\text{Sect}(M)$ is the quotient of the semigroup of the
endomorphisms of $M$ modulo the equivalence relation: $\rho,\rho'\in
\text{End}(M),\, \rho\thicksim\rho'$ iff there is a unitary $u\in M$
such that $\rho'(x)=u\rho(x)u^*$ for all $x\in M$.

$\text{Sect}(M)$ is a $^*$-semiring (there are an addition, a
product and an involution $\rho\rightarrow \bar\rho$) equivalent to
the Connes correspondences (bimodules) on $M$ up to unitary
equivalence. If $\r$ is an element of $\text{End}(M)$ we shall
denote by $[\r]$ its class in $\text{Sect}(M)$. We define
$\text{Hom}(\r,\r')$ between the objects $\r,\r'\in \End(M)$ by
\[
\text{Hom}(\r,\r')\equiv\{a\in M: a\r(x)=\r'(x)a \ \forall x\in M\}.
\]
We use $\langle  \lambda , \mu \rangle$ to denote the dimension of
$\text{\rm Hom}(\lambda , \mu )$; it can be $\infty$, but it is
finite if $\la,\mu$ have finite index. See \cite{J1} for the
definition of index for type $II_1$ case which initiated the subject
and  \cite{PP} for  the definition of index in general. Also see
\S2.3 \cite{KLX} for expositions. $\langle  \lambda , \mu \rangle$
depends only on $[\lambda ]$ and $[\mu ]$. Moreover we have if $\nu$
has finite index, then $\langle \nu \lambda , \mu \rangle = \langle
\lambda , \bar \nu \mu \rangle $, $\langle \lambda\nu , \mu \rangle
= \langle \lambda , \mu \bar \nu \rangle $ which follows from
Frobenius duality. $\mu $ is a subsector of $\lambda $ if there is
an isometry $v\in M$ such that $\mu(x)= v^* \lambda(x)v, \forall
x\in M.$ We will also use the following notation: if $\mu $ is a
subsector of $\lambda $, we will write as $\mu \prec \lambda $  or
$\lambda \succ \mu $.  A sector is said to be irreducible if it has
only one subsector.

\subsection{Local nets}
%\subsubsection{Preliminaries on conformal nets}
By an interval of the circle we mean an open connected non-empty
subset $I$ of $S^1$ such that the interior of its complement $I'$ is
not empty. We denote by $\I$ the family of all intervals of $S^1$.

A {\it net} $\A$ of von Neumann algebras on $S^1$ is a map
\[
I\in\I\to\A(I)\subset B(\H)
\]
from $\I$ to von Neumann algebras on a fixed separable Hilbert space
$\H$ that satisfies:
\begin{itemize}
\item[{\bf A.}] {\it Isotony}. If $I_{1}\subset I_{2}$ belong to
$\I$, then
\begin{equation*}
 \A(I_{1})\subset\A(I_{2}).
\end{equation*}
\end{itemize}
If $E\subset S^1$ is any region, we shall put
$\A(E)\equiv\bigvee_{E\supset I\in\I}\A(I)$ with $\A(E)=\mathbb C$
if $E$ has empty interior (the symbol $\vee$ denotes the von Neumann
algebra generated).

The net $\A$ is called {\it local} if it satisfies:
\begin{itemize}
\item[{\bf B.}] {\it Locality}. If $I_{1},I_{2}\in\I$ and $I_1\cap
I_2=\emptyset$ then
\begin{equation*}
 [\A(I_{1}),\A(I_{2})]=\{0\},
 \end{equation*}
where brackets denote the commutator.
\end{itemize}
The net $\A$ is called {\it M\"{o}bius covariant} if in addition
satisfies the following properties {\bf C,D,E,F}:
\begin{itemize}
\item[{\bf C.}] {\it M\"{o}bius covariance}.
There exists a non-trivial strongly continuous unitary
representation $U$ of the M\"{o}bius group $\Mob$ (isomorphic to
$PSU(1,1)$) on $\H$ such that
\begin{equation*}
 U(g)\A(I) U(g)^*\ =\ \A(gI),\quad g\in \Mob,\ I\in\I.
\end{equation*}
\item[{\bf D.}] {\it Positivity of the energy}.
The generator of the one-parameter rotation subgroup of $U$
(conformal Hamiltonian), denoted by $L_0$ in the following,  is
positive.
\item[{\bf E.}] {\it Existence of the vacuum}.  There exists a unit
$U$-invariant vector $\Omega\in\H$ (vacuum vector), and $\Omega$ is
cyclic for the von Neumann algebra $\bigvee_{I\in\I}\A(I)$.
\end{itemize}
By the Reeh-Schlieder theorem $\Omega$ is cyclic and separating for
every fixed $\A(I)$. The modular objects associated with
$(\A(I),\Omega)$ have a geometric meaning
\[
\Delta^{it}_I = U(\Lambda_I(2\pi t)),\qquad J_I = U(r_I)\ .
\]
Here $\Lambda_I$ is a canonical one-parameter subgroup of $\Mob$ and
$U(r_I)$ is a antiunitary acting geometrically on $\A$ as a
reflection $r_I$ on $S^1$.

This implies {\em Haag duality}:
\[
\A(I)'=\A(I'),\quad I\in\I\ ,
\]
where $I'$ is the interior of $S^1\setminus I$.

\begin{itemize}
\item[{\bf F.}] {\it Irreducibility}. $\bigvee_{I\in\I}\A(I)=B(\H)$.
Indeed $\A$ is irreducible iff $\Om$ is the unique $U$-invariant
vector (up to scalar multiples). Also  $\A$ is irreducible iff the
local von Neumann algebras $\A(I)$ are factors. In this case they
are either ${\mathbb C}$ or III$_1$-factors with separable predual
in Connes classification of type III factors.
\end{itemize}
By a {\it conformal net} (or diffeomorphism covariant net) $\A$ we
shall mean a M\"{o}bius covariant net such that the following holds:
\begin{itemize}
\item[{\bf G.}] {\it Conformal covariance}. There exists a projective
unitary representation $U$ of $\Diff(S^1)$ on $\H$ extending the
unitary representation of $\Mob$ such that for all $I\in\I$ we have
\begin{gather*}
 U(\phi)\A(I) U(\phi)^*\ =\ \A(\phi.I),\quad  \phi\in\Diff(S^1), \\
 U(\phi)xU(\phi)^*\ =\ x,\quad x\in\A(I),\ \phi\in\Diff(I'),
\end{gather*}
\end{itemize}
where $\Diff(S^1)$ denotes the group of smooth, positively oriented
diffeomorphism of $S^1$ and $\Diff(I)$ the subgroup of
diffeomorphisms $g$ such that $\phi(z)=z$ for all $z\in I'$.
\par
A (DHR) representation $\pi$ of $\A$ on a Hilbert space $\H$ is a
map $I\in\I\mapsto  \pi_I$ that associates to each $I$ a normal
representation of $\A(I)$ on $B(\H)$ such that
\[
\pi_{\tilde I}\res\A(I)=\pi_I,\quad I\subset\tilde I, \quad I,\tilde
I\subset\I\ .
\]
$\pi$ is said to be M\"obius (resp. diffeomorphism) covariant if
there is a projective unitary representation $U_{\pi}$ of $\Mob$
(resp. $\Diff(S^1)$) on $\H$ such that
\[
\pi_{gI}(U(g)xU(g)^*) =U_{\pi}(g)\pi_{I}(x)U_{\pi}(g)^*
\]
for all $I\in\I$, $x\in\A(I)$ and $g\in \Mob$ (resp.
$g\in\Diff(S^1)$).

By definition the irreducible conformal net is in fact an
irreducible representation of itself and we will call this
representation the {\it vacuum representation}.\par

Let $G$ be a simply connected  compact Lie group. By Th. 3.2 of
\cite{FG}, the vacuum positive energy representation of the loop
group $LG$ (cf. \cite{PS}) at level $k$ gives rise to an irreducible
conformal net denoted by {\it ${\A}_{G_k}$}. By Th. 3.3 of
\cite{FG}, every irreducible positive energy representation of the
loop group $LG$ at level $k$ gives rise to  an irreducible covariant
representation of ${\A}_{G_k}$. \par Given an interval $I$ and a
representation $\pi$ of $\A$, there is an {\em endomorphism of $\A$
localized in $I$} equivalent to $\pi$; namely $\r$ is a
representation of $\A$ on the vacuum Hilbert space $\H$, unitarily
equivalent to $\pi$, such that
$\r_{I'}=\text{id}\restriction\A(I')$. We now define  the
statistics. Given the endomorphism $\r$ of $\A$ localized in
$I\in\I$, choose an equivalent endomorphism $\r_0$ localized in an
interval $I_0\in\I$ with $\bar I_0\cap\bar I =\emptyset$ and let $u$
be a local intertwiner in $\Hom(\r,\r_0)$ , namely $u\in
\Hom(\r_{\tilde I},\r_{0,\tilde I})$ with $I_0$ following clockwise
$I$ inside $\tilde I$ which is an interval containing both $I$ and
$I_0$.

The {\it statistics operator} $\epsilon (\r,\rho):= u^*\r(u) =
u^*\r_{\tilde I}(u) $ belongs to $\Hom(\r^2_{\tilde I},\r^2_{\tilde
I})$. We will call $\epsilon (\r,\rho)$ the positive or right
braiding and $\tilde\epsilon (\r,\rho):=\epsilon (\r,\rho)^*$ the
negative or left braiding.
\par
Let $\B$ be a  conformal  net. By a {\it  conformal subnet} (cf.
\cite{Longo4}) we shall mean a map
\[
I\in\I\to\A(I)\subset \B(I)
\]
that associates to each interval $I\in \I$ a von Neumann subalgebra
$\A(I)$ of $\B(I)$, which is isotonic
\[
\A(I_1)\subset \A(I_2), I_1\subset I_2,
\]
and  conformal covariant with respect to the representation $U$,
namely
\[
U(g) \A(I) U(g)^*= \A(g.I)
\] for all $g\in \Diff(S^1)$ and $I\in \I$. Note that by Lemma 13
of \cite{Longo4} for each $I\in \I$ there exists a conditional
expectation $E_I: \B(I)\rightarrow \A(I)$ such that $E_I$ preserves
the vector state given by the vacuum of $\A$.
\begin{definition}\label{ext}
Let $\A$ be a  conformal net. A  conformal net $\B$ on a Hilbert
space $\H$ is an extension of $\A$ or $\A$ is a subnet of $\B$ if
there is a DHR representation $\pi$ of $\A$ on $\H$ such that
$\pi(\A)\subset \B$ is a conformal subnet. The extension is
irreducible if $\pi(\A(I))'\cap \B(I) = {\mathbb C} $ for some (and
hence all) interval $I$, and is of finite index if
$\pi(\A(I))\subset \B(I)$ has finite index for some (and hence all)
interval $I$. The index will be called the index of the inclusion
$\pi(\A)\subset \B$ and is denoted by $[\B:\A].$ If $\pi$ as
representation of
 $\A$ decomposes as $[\pi]= \sum_\lambda m_\lambda[\lambda]$ where
$m_\lambda$  are non-negative  integers and $\lambda$ are
irreducible DHR representations of $\A$, we say that $[\pi]=
\sum_\lambda m_\lambda[\lambda]$ is the spectrum of the extension.
For simplicity we will write  $\pi(\A)\subset \B$ simply as
$\A\subset \B$.
\end{definition}
\begin{lemma}
If $\A\subset \B$ is a conformal subnet with finite index, then
$\A\subset \B$ is irreducible.
\end{lemma}
\proof This is proved in Cor. 3.6 of \cite{BE1}, without assumption
of conformal covariance of $\A$ but under the additional assumption
that $\A$ is strongly additive to ensure the equivalence of local
and global intertwiners, but for conformal net $\A$ the equivalence
of local and global intertwiners for finite index representations
are proved in \S2 of \cite{GL2}, thus the proof of Cor. 3.6 of
\cite{BE1} applies verbatim.
\endproof

\begin{lemma}\label{charge}
Suppose that $\A\subset \B$ has finite index , and let $[\pi]=
\sum_\lambda m_\lambda[\lambda]$ be as in Definition above. Fix an
interval $I$ and suppose that $\lambda,\bar\la$ is localized on $I.$
\par

(1) Let $K_\lambda:=\{ T\in B(I)| Ta=a\lambda(a)T, \forall a\in
\A(I) \}.$ Then $K_\la$ is a vector space of dimension $m_\la\leq
d_\la.$  One can find isometries $T_{\lambda_i}\in K_\lambda,
T_{\lambda_i}\in K_{\bar\lambda}, 1\leq i\leq m_\lambda$ such that
$T_{\lambda_i} a= \lambda(a) T_{\lambda_i}, \forall a\in \A,
E(T_{\lambda_i}T_{\lambda_j^*})= \delta_{ij} 1/d_\lambda,
E(T_{\bar\lambda_i}T_{\bar\lambda_j^*})= \delta_{ij}
1/d_\lambda,T_{\la_i}^*\in \A(I)T_{\bar\la_i};$ Every $b\in \B(I)$
can be written as $b=\sum_{\lambda_i}
{d_\la}T_{\lambda_i}^*E(T_{\lambda_i}b);$
\par
(2) Let $L_\lambda \subset K_\la$ be subspaces with the following
properties:(a) $L_\lambda L_\mu \subset \sum_\nu \A(I) L_\nu;$ (b)
$L_\la^*\subset \A(I) L_{\bar\la}.$ Then there is an intermediate
subnet $\A\subset \C\subset \B$ such that $\C(I)=\sum_\lambda
\A(I)L_\la.$ Conversely every intermediate subnet arises this way;

\par

(3) If $\Omega$ is the vacuum vector of $\B,$ and denote by
$\overline{\A\Omega}=H_0, \overline{
{T_{\lambda_i}^*}\A\Omega}=H_{\lambda_i},$ then as Hilbert space
$H=\bigoplus_{\lambda_i,1\leq i\leq m_\lambda} H_{\lambda_i},$ and
the map $\sqrt{d_\lambda}T_{\lambda_i}^*: H_0\rightarrow
H_{\lambda_i}$ is a unitary intertwiner between the action of
$\lambda(\A(I))$ on $H_0$ and $\A(I)$ on $H_{\lambda_i}.$

\end{lemma}
\proof

(1) and(2) follow from \S3 of \cite{ILP} and \S2 of \cite{Longo4}.
For (3), only unitarity has to be checked. We have
$$
\langle T_{\lambda_i}^* a_1\Omega, T_{\lambda_i}^* a_2\Omega\rangle
=\langle a_2^*E(T_{\lambda_i}T_{\lambda_i}^*)a_1\Omega,
\Omega\rangle= 1/d_{\lambda}\langle a_1\Omega,a_2\Omega\rangle,
\forall a_1,a_2\in \A(I),
$$
and the proof is complete.
\endproof

\subsection{Induced endomorphisms}\label{inductionsection}
 Suppose a conformal net
$\A$ and a representation $\lambda$ is given. Fix an open interval I
of the circle and Let $M:=\A(I)$ be a fixed type $III_1$ factor.
Then $\lambda$ give rises to an endomorphism still denoted by
$\lambda$ of $M$. Suppose $\{[\lambda] \}$ is a finite set  of all
equivalence classes of irreducible, covariant, finite-index
representations of an irreducible local conformal net $\A$. We will
use $\Delta_\A$ to denote all finite index representations of net
$\A$ and will use the same notation $\Delta_\A$ to denote the
corresponding sectors of $M$.

We will denote the conjugate of $[\lambda]$ by $[{\bar \lambda}]$
and identity sector (corresponding to the vacuum representation) by
$[1]$ if no confusion arises, and let $N_{\lambda\mu}^\nu = \langle
[\lambda][\mu], [\nu]\rangle $. Here $\langle \mu,\nu\rangle$
denotes the dimension of the space of intertwiners from $\mu$ to
$\nu$ (denoted by $\text {\rm Hom}(\mu,\nu)$).  The univalence of
$\lambda$ and the statistical dimension of (cf. \S2  of \cite{GL2})
will be denoted by $\omega_{\lambda}$ and $d{(\lambda)}$ (or
$d_{\lambda})$) respectively. Suppose that $\ro\in \End(M)$ has the
property that $\gamma=\rho\bar\rho \in \Delta_\A$. By \S2.7 of
\cite{LR}, we can find two isometries $v_1\in \Hom(\gamma,\gamma^2),
w_1\in \Hom (1,\gamma)$ such that $\bar\rho(M)$ and $v_1$ generate
$M$ and
\begin{align*}
v_1^* w_1 & = v_1^* \gamma(w_1) = d_\rho^{-1} \\
v_1v_1 & = \gamma(v_1) v_1  \\
%\epsilon(\ra,\r) w_1 & = w_1 \label{c}
\end{align*}
By Thm. 4.9 of \cite{LR}, we shall say that $\rho$ is {\it local }
if
\begin{align}\label{local}
v_1^* w_1 & = v_1^* \gamma(w_1) = d_\rho^{-1} \\
v_1v_1 & = \gamma(v_1) v_1  \\
\bar\rho(\epsilon(\gamma,\gamma)) v_1 & = v_1
\end{align}
Note that if $\rho$ is local, then
\begin{equation}\label{local=1}
\om_\mu =1, \forall \mu\prec \ro\bar\ro \end{equation}

For each (not necessarily irreducible)  $\lambda\in \Delta_\A,$ let
$\e{(\lambda,\gamma)}
 :   \lambda\gamma \rightarrow
\gamma\lambda$ (resp.  $\tilde\e{(\lambda,\gamma)}$), be the
positive (resp. negative) braiding operator as defined in Section
1.4 of \cite{Xb}. Denote by $\lambda_\e \in$ End$(M)$ which is
defined by
\begin{align*}
\lambda_\e(x) :&= ad(\e(\lambda,\gamma))\lambda(x)=
\e(\lambda,\gamma) \lambda(x) \e(\lambda,\gamma)^* \\
\lambda_{\tilde\e}(x) :&= ad(\tilde\e(\lambda,\gamma))\lambda(x)=
\tilde\e(\lambda,\gamma)^* \lambda(x) \tilde\e(\lambda,\gamma)^*,
\forall x\in M.
\end{align*}
 By (1) of Theorem 3.1 of \cite{Xb},  $\lambda_{\e} \rho (M)
\subset \rho(M), \lambda_{\tilde\e}  \rho (M) \subset \rho(M),$
hence the following definition makes sense:
 \bdef\label{ala} If
$\la\in \Delta_\A$ define two elements of $\End(M)$ by
$$
a^\ro_\lambda(m):= \rho^{-1} (\lambda_\e  \rho (m)) , \ \tilde
a^\ro_\lambda(m):= \rho^{-1} (\lambda_{\tilde\e}  \rho (m)), \forall
m\in M.
$$
$a_\lambda^\ro$ (resp. $\tilde a_\lambda^\ro$) will be referred to
as positive (resp. negative) induction of $\la$ with respect to
$\ro.$

\ede
\begin{remark}
For simplicity we will use $a_\la, \tilde a_\la$ to denote
$a_\lambda^\ro, \tilde a_\lambda^\ro$ when it is clear that
inductions are with respect to the same $\ro.$
\end{remark}

The endomorphisms $a_\lambda$ are called braided endomorphisms in
\cite{Xb}  due to its braiding properties (cf.  (2) of Corollary 3.4
in \cite{Xb}), and enjoy an interesting set of properties (cf.
Section 3 of \cite{Xb}).  We summarize a few properties from
\cite{Xb} which will be used in this paper: (cf. Th. 3.1 , Co. 3.2
and Th. 3.3 of \cite{Xb} ): \bprop\label{xua} (1). The maps
$[\lambda] \rightarrow [a_\lambda], [\lambda]\rightarrow [\tilde
a_\lambda]$ are ring homomorphisms;\par (2) $ a_\lambda \bar \rho
=\tilde a_\lambda \bar \rho =\bar \rho \lambda$;\par (3) When
$\ro\bar\ro$ is local, $ \langle a_\la,  a_\mu \ra = \langle \tilde
a_\la, \tilde a_\mu \ra = \langle a_\la\bar \ro, a_\mu\bar\ro \ra
=\langle \tilde a_\la\bar \ro, \tilde a_\mu\bar\ro \ra ;$\par (4)
(3) remains valid if $a_\la, a_\mu$ (resp. $\tilde a_\la, \tilde
a_\mu$) are replaced by their subsectors. In particular we have
$\langle a_\la, \sigma \rangle = \langle \la,\rho\sigma\bar\rho
\rangle$  if $\sigma\prec a_\mu.$

\eprop

The following is Porp. 2.24 of \cite{X-adv}:
 \bprop\label{loclr}
Suppose that $\ro\bar\ro \in \Delta$. Then:\par (1) $\ro$ is local
iff $\langle 1, a_\mu\ra =\langle \ro \bar\ro, \mu\ra, \forall
\mu\in \Delta_\A;$
\par (2) \[ \ro= \ro'\ro''= \tilde\ro'\tilde\ro'' \] where
$\ro',\ro'', \tilde\ro', \tilde\ro'' \in \End (M)$, and
$\ro',\tilde\ro'$ are local which verifies
\begin{align*}
\langle \ro'\bar\ro',\mu\ra &= \langle 1,a_\mu\ra =\langle 1,
a_\mu^{\ro'}\ra
\\
\langle \tilde\ro'\overline{\tilde\ro'},\mu\ra &= \langle 1,\tilde
a_\mu\ra =\langle 1, \tilde a_\mu^{\tilde \ro'}\ra
\end{align*}
$\forall \mu\in \Delta_\A.$ We refer to $\rho'$ (resp. $\rho''$) as
the left (resp.right) {\it local support} of  $\rho$.

\eprop The following Lemma is Prop. 3.23 of \cite{BE1} (The proof
was also implicitly contained in the proof of Lemma 3.2 of
\cite{Xb}):

\blem\label{loc2} If $\rho\bar\ro$ is local, then $[a_\la]=[\tilde
a_\la]$ iff $\e(\la,\rho\bar\rho)\e(\ro\bar\ro,\la) =1$ iff
$\e(\la,\mu)\e(\mu,\la) =1, \forall \mu\in \ro\bar\rho.$ \elem We
shall make use of the following notation in \S4:
\bdef\label{zmatrix} For $\la,\mu\in \Delta_\A$, $Z_{\la\mu}^\rho:=
\langle a_\lambda, \tilde a_\mu\ra.$ \ede
%\subsection{Type A}
\subsection {Jones-Wassermann subfactors from representation of Loop
groups}\label{typea} Let $G= SU(n)$. We denote $LG$ the group of
smooth maps $f: S^1 \mapsto G$ under pointwise multiplication. The
diffeomorphism group of the circle $\text{\rm Diff} S^1 $ is
naturally a subgroup of $\text{\rm Aut}(LG)$ with the action given
by reparametrization. In particular the group of rotations
$\text{\rm Rot}S^1 \simeq U(1)$ acts on $LG$. We will be interested
in the projective unitary representation $\pi : LG \rightarrow U(H)$
that are both irreducible and have positive energy. This means that
$\pi $ should extend to $LG\rtimes \text{\rm Rot}\ S^1$ so that
$H=\oplus _{n\geq 0} H(n)$, where the $H(n)$ are the eigenspace for
the action of $\text{\rm Rot}S^1$, i.e., $r_\theta \xi = \exp(i n
\theta)$ for $\theta \in H(n)$ and $\text{\rm dim}\ H(n) < \infty $
with $H(0) \neq 0$. It follows from \cite{PS} that for fixed level
$k$ which is a positive integer, there are only finite number of
such irreducible representations indexed by the finite set
$$
 P_{++}^{k}
= \bigg \{ \lambda \in P \mid \lambda = \sum _{i=1, \cdots , n-1}
\lambda _i \Lambda _i , \lambda _i \geq 0\, , \sum _{i=1, \cdots ,
n-1} \lambda _i \leq k \bigg \}
$$
where $P$ is the weight lattice of $SU(n)$ and $\Lambda _i$ are the
fundamental weights. We will write $\la=(\la_1,...,\la_{n-1}),
\la_0= k-\sum_{1\leq i\leq n-1} \la_i$ and refer to
$\la_0,...,\la_{n-1}$ as components of $\la.$

We will use $\Lambda_0$ or simply $1$  to denote the trivial
representation of $SU(n)$. For $\lambda , \mu , \nu \in P_{++}^{k}$,
define $N_{\lambda \mu}^\nu  = \sum _{\delta \in P_{++}^{k}
}S_\lambda ^{(\delta)} S_\mu ^{(\delta)} S_\nu
^{(\delta*)}/S_{\Lambda_0}^{(\delta})$ where $S_\lambda ^{(\delta)}$
is given by the Kac-Peterson formula:
$$
S_\lambda ^{(\delta)} = c \sum _{w\in S_n} \varepsilon _w \exp
(iw(\delta) \cdot \lambda 2 \pi /n)
$$
where $\varepsilon _w = \text{\rm det}(w)$ and $c$ is a
normalization constant fixed by the requirement that
$S_\mu^{(\delta)}$ is an orthonormal system. It is shown in
\cite{Kac2} P. 288 that $N_{\lambda \mu}^\nu $ are non-negative
integers. Moreover, define $ Gr(C_k)$ to be the ring whose basis are
elements of $ P_{++}^{k}$ with structure constants $N_{\lambda
\mu}^\nu $.
  The natural involution $*$ on $ P_{++}^{k}$ is
defined by $\lambda \mapsto \lambda ^* =$ the conjugate of $\lambda
$ as representation of $SU(n)$. Note that $\lambda\rightarrow
\frac{S_{\la\mu}}{S_{1\mu}}$ gives a representation of $ Gr(C_k).$

\par

We shall also denote $S_{\Lambda _0}^{(\Lambda)}$ by $S_1^{(\Lambda
)}$. Define $d_\lambda = \frac {S_1^{(\lambda )}}{S_1^{(\Lambda
_0)}}$. We shall call $(S_\nu ^{(\delta )})$ the $S$-matrix of
$LSU(n)$ at level $k$. \par
 We shall encounter the $\Bbb Z_n$
group of automorphisms of this set of weights, generated by
$$
\sigma : \lambda = (\lambda_1, \lambda_2, \cdots , \lambda_{n-1})
\rightarrow \sigma(\lambda) = ( k -1- \lambda_1 -\cdots
\lambda_{n-1}, \lambda_1, \cdots , \lambda_{n-2}).
$$
Define  $\col(\lambda) = \Sigma_i (\lambda_i - 1) i $. $\col(\la)$
will be referred to as the color of $\la$.  The central element  $
\exp \frac{2\pi i}{n}$ of $SU(n)$ acts on representation of $SU(n)$
labeled by $\lambda$ as $\exp( \frac{2\pi i \col(\lambda)}{n})$. The
irreducible positive energy representations of $ L SU(n)$ at level
$k$ give rise to an irreducible conformal net $\A$ (cf. \cite{KLX})
and its covariant representations. We will use
$\la=(\la_1,...\la_{n-1})$ to denote irreducible representations of
$\A$ and also the corresponding endomorphism of $M=\A(I).$

All the sectors $[\lambda]$ with $\lambda$ irreducible generate the
fusion ring of $\A.$
\par For $\lambda$ irreducible, the univalence $\om_\lambda$ is given
by an explicit formula (cf. 9.4 of [PS]). Let us first define
$h_\lambda = \frac {c_2(\lambda)}{k+n}$ where $c_2(\lambda)$ is the
value of Casimir operator on representation of $SU(n)$ labeled by
dominant weight $\lambda$.
 $h_\lambda$ is usually called the conformal dimension. Then
we have: $\om_\lambda = exp({2\pi i} h_\lambda)$. The conformal
dimension of $\lambda=(\la_1,...,\la_{n-1})$ is given by
\begin{equation}\label{cdim} h_\lambda= \frac{1}{2n(k+n)}\sum_{1\leq
i\leq n-1} i(n-i) \la_i^2 + \frac{1}{n(k+n)}\sum_{1\leq j\leq i\leq
n-1}j (n-i)\la_j\la_i + \frac{1}{2(k+n)}\sum_{1\leq j\leq n-1}
j(n-j) \la_j \end{equation} The following  result is proved in
\cite{Wass} (See Corollary 1 of Chapter V in \cite{Wass}).

\bthm\label{wass}  Each $\lambda \in  P_{++}^{(k)}$ has finite index
with index value $d_\lambda ^2$.  The fusion ring generated by all
$\lambda \in P_{++}^{(k)}$ is isomorphic to $ Gr(C_k)$. \ethm

\begin{remark}\label{jxx}
The subfactors in the above theorem are called Jones-Wassermann
subfactors after the authors who first studied them (cf.
\cite{Jh1},\cite{Wass}).
\end{remark}
\begin{definition}
$v:=(1,0,...,0), \ad:=(1,0,...,0,1), \om^i =k \Lambda_i, 0\leq i\leq
n-1. $  $v$ (resp. $v_0$) will be referred to as vector (resp.
adjoint) representation. \ede The following is observed in
\cite{GW}: \blem\label{fusionrule} Let $(0,...,0,1,0,...0)$ be the
$i$-th ($1\leq i\leq n-1$) fundamental weight. Then
$[(0,...,0,1,0,...0)\la]$ are determined as follows:
$\mu\prec(0,...,0,1,0,...0)\la$ iff when the Young diagram of $\mu$
can be obtained from Young diagram of $\la$ by adding $i$ boxes on
$i$ different rows of $\la$, and such $\mu$ appears in
$[(0,...,0,1,0,...0)\la]$ only once. \elem

 \blem\label{fusionv}
(1) If $[\la]\neq \omega^i$ for some $0\leq i\leq n-1,$ then
$\ad\prec \la\bar\la;$\par (2) If $\lambda_1\lambda_2$ is
irreducible, then either $\lambda_1$ or $\lambda_2=\omega^i$ for
some $0\leq i\leq n-1;$\par (3) Suppose that $\la$ has color $0\mod
n.$ Then $\la\prec v_0^m$ for some $m \in \mathbb{N}.$ \elem \proof
(1), (2) is lemma 2.30 of \cite{X-adv}. By the lemma above $\la\prec
v^l$ for some $l\in \mathbb{N},$ and since $\col(\la)=0 \ \mod  n,$
we have $l=nl_1, l_1\in \mathbb{N}.$ Since $1\prec v^n, 1\prec {\bar
v}^n,$ we have $[v^n]\prec [v^n\bar v^n] =([v_0]+[1])^n,$ and (3)
follows.
\endproof
\subsection{Subnets from conformal inclusions}\label{ci}

 Let $G \subset H$ be inclusions of compact
simply connected Lie groups.  $LG \subset LH$ is called a conformal
inclusion if the level 1 projective positive energy representations
of $LH$ decompose as a finite number of irreducible projective
representations of $LG$.   $LG \subset LH$ is called a maximal
conformal inclusion if there is no proper subgroup $G'$ of $H$
containing $G$ such that   $LG \subset LG'$ is also  a conformal
inclusion. A list of maximal conformal inclusions can be found in
\cite{GNO}.  \par Let $H^0$ be the vacuum representation of $LH$,
i.e., the
 representation of $LH$ associated with the trivial  representation of $H$.
Then $H^0$ decomposes as a direct sum of  irreducible projective
representation of $LG$ at level $K$. $K$  is called the Dynkin index
of the conformal inclusion.

%Assume $H^0 = \oplus_{\lambda \in P_0} m_\lambda
%H_\lambda$ where $P_0 \subset C_K$ is finite and $m_\lambda$ is the
%multiplicity of $H_\lambda$ in $H^0$ .
We shall  write the conformal inclusion as $G_K\subset H_1$. Note
that it follows from the definition that $\A_{H_1}$ is an extension
of $\A_{G_K}$. We shall limit our consideration to the following
conformal inclusions so we can use the results of \cite{Xb}:
\begin{align}
{SU}(n)_{n-2} & \subset \ {SU} \left( \frac{n(n-1)}{2}
\right)_1, \ \ N \geq 4 ; \label{c1} \\
{SU}(n)_{n+2} & \subset \ {SU} \left( \frac{n(n+1)}{2}
\right)_1;  \label{c2} \\
{SU}(n)_{n} & \subset \ {Spin} (n^2 - 1)_1, \ \ N \geq
2;  \label{c3} \\
{SU}(n)_m \times {SU}(m)_n  & \subset \ {SU(mn)}_1. \label{c4}
\end{align}
Note that except equation (\ref{c4}), the above cover all the
infinite series  of maximal conformal inclusions of the form $SU(N)
\subset H$ with $H$ being a simple group.

\section{Intermediate subnets in confonmal subnets associated with
conformal inclusions}

Let $\A\subset \B$ be conformal subnets associated with conformal
inclusions in \S\ref{ci}, i.e., $\A=\A_{G_k}\subset \B=\A_{H_1}.$
Our goal in this section is to list all intermediate subnets
$\A\subset \C\subset \B.$ \par The spectrum $[\pi]=\sum_\lambda
m_\lambda \lambda$ of $\A\subset \B$ is given by \cite{ABI} and
\cite{LL}. One interesting feature is that all $m_\lambda=1.$ We
write $H_\B= \oplus_\lambda H_\lambda$ with $H_0$ the vacuum
representation of $\A,$ and $H_\B$ (resp. $H_\C$) the vacuum
representation space of $\B$ (resp. $H_\C$).\par Fix an interval $I$
and let $M=\A(I)\subset \C(I), \rho\in \End(M), \rho\bar\rho=H_\C\in
\Delta_\A$ where we use $H_\C$ to denote the restriction of the
vacuum representation of $\C$ to $\A.$ For $\la\in \Delta_\A,$ we
will write $a_\la^\C:=a_\la^\rho.$
\subsection{Smeared Vertex Operators}\label{smear}
Let $g$ (resp. $h$) be the Lie algebra of $G$ (resp. $H$). Choose a
basis $e_\alpha,e_{-\alpha},h_\alpha$ in $h_{\Bbb C}:=h\otimes {\Bbb
C}$ with $\alpha$ ranging over the set of roots as in \S 2.5 of
\cite{PS}. Let $X_\alpha :=i(e_\alpha + e_{-\alpha}), Y_\alpha:=
(e_\alpha -e_{-\alpha})$. Denote by $\hat h$ the affine Kac-Moody
algebra (cf. P. 163 of \cite{[KW]}) associated to $h_{\Bbb C}$. Note
$\hat h=h_{\Bbb C}\otimes {\Bbb C}[t,t^{-1}] \oplus {\Bbb C} c$,
where ${\Bbb C}c$ is the 1-dimensional center of $\hat h$. For $X\in
h$, Define $X(n):= X \otimes t^n$,  $X(z):= \sum_n X(n) z^{-n-1}$ as
on Page 312 of \cite{[KT]}.
\par
Let $\pi^0$ be the vacuum representation of $LH_1$ on $H_\B$ with
vacuum vector $\Omega$.  Let $D$ be the generator of the action of
the rotation group on $H_\B$.  $H^0_\B$ will denote the finite
linear sum of the eigenvectors of $D$. For $\xi\in H_\B$, we define
$||x||_s = ||(1+D)^s x||, s\in {\Bbb R}$. $H^{s}:=\{x\in H^0_\B |\
||x||_s < \infty \}$ and $H(\infty) = \cap_{s\in {\Bbb R}} H^{s}$.
Note that when $s\geq 0$, $H^{s}$ is a complete space under the norm
$||.||_s$.  Clearly $H^0_\B \subset H(\infty)$. The elements of
$H^0_\B$ (resp. $H(\infty)$) will be called {\it finite energy
vectors} (resp. {\it smooth vectors}). The eigenvalue of $D$ is
sometimes referred to as energy or weight. \par Let us recall a few
elementary facts about vertex operators which will be used. See
\cite{[FLM]} or \cite{[Kacv]}  for an introduction on vertex
operator algebras. Define $End(H^0_\B)$ to be the space of all
linear operators (not necessarily bounded) from $H^0_\B$ to $H^0_\B$
and set
$$
End(H^0_\B)[[z,z^{-1}]]:= \{ \sum_{n\in {\Bbb Z}} v_n z^n| v_n\in
End(H^0_\B) \}.
$$
By the statement on P. 154 of \cite{[FZ]} which follows from Th.
2.4.1 of \cite{[FZ]} there exists a linear map
$$
\psi\in H^0_\B \rightarrow V(\xi,z)= \sum_{m\in {\Bbb Z}} \psi(m)
z^{-m-1} \in End(H^0_\B)[[z,z^{-1}]]
$$
with the following properties: \par (1) $\psi(-1) \Omega = \psi$;
\par (2) If
$$
\psi= X_{i_1}(-1)... X_{i_t}(-1) \Omega ,$$ then
$$
V(\psi,z)=  :X_{i_1}(z)... X_{i_t}(z):
$$ where $:,:$ are normal ordered products (cf. (2.38), (2.39) of \cite{D}).
\par
$V(\psi,z)$ is called a {\it vertex operator} of $\psi$. \par

%Recall $ V(\psi, z)= \sum_m \psi(m) z^{-m-1}$. Define
%$$
%V(m):= \psi(m+n-1)
%$$
%so we have $ V(\psi, z)= \sum_m V(m)  z^{-m-n}$. This  expression
%for $ V(\psi, z)$ is in accordance with the convention of [KT].
%Note that $V(-n)\Omega = \psi(-1)\Omega = \psi$ by property (1)
%above.
\par
Let $f=\sum_m f(m) z^m$ be a smooth test function. Define
$$
||f||_s= \sum_{n\in {\Bbb Z}} (1+|m|)^s |f(m)|.
$$
The {\it smeared vertex operator} $V(\psi,f)$ is defined to be:
$$
V(\psi,f) = \frac{1}{2\pi i}\int_{S^1} V( \psi, z) f dz = \sum_m
f(m) \psi(m) .$$ $V(\psi,f)$ is a well defined operator on $H^0_\B$.
Let $V(\psi,f)^{FA}$ be the formal adjoint of $V(\psi,f)$ on
$H^0_\B$. It is defined by the equation
$$
\langle V(\psi,f)x, y\rangle = \langle x, V(\psi,f)^{FA}y \rangle ,
\forall x,y\in H^0_\B
$$
where $\langle, \rangle$ is the inner product on Hilbert space
$H^0_\B$.
%Similarly for $X\in h$, we define $X(f):=\sum_n
%X(n)f(n)$.\par

\begin{lemma}\label{RS}
The subspace spanned by $V(\psi,f)\Omega, \forall \psi\in
H_\lambda^0=H_\la\cap H_\B^0, \forall f$ smooth, $supp f\in I,$  is
dense in $H_\lambda.$
\end{lemma}
\proof The proof is essentially the same as the proof of
Reeh-Schlieder Th. Let $\xi\in H_\lambda$ be a vector which is
orthogonal to the subspace spanned by $V(\psi,f)\Omega, \forall
\psi\in H_\lambda^0, supp f\in I.$ Suppose $J$ is an open interval
such that $\bar J\subset I,$ and $f$ is a smooth function with
support in $J$. Consider the function
$$
F(z)=\langle \exp(izD)\xi, V(\psi,f)\Omega\rangle.
$$
Since the spectrum of $D$ on $H_\lambda$ is a subset of non-negative
integers, it follows that $F(z)$ is holomorphic on the upper half
plane, continues on the real line, and vanishes on an open interval
on the real line. It follows by Schwartz reflection principle that
$F(z)$ is identically zero, and we have
$$
\langle \exp(itD)\xi, V(\psi,f)\Omega\rangle=\langle \xi,
\exp(-itD)V(\psi,f)\exp(it D) \Omega\rangle=0, \forall t\in
\mathbb{R}.
$$
On the other hand $\exp(-itD)V(\psi,f)\exp(it D) \Omega= V(\psi,R_t
(f))\Omega,$  where $R_t(f)(z)= f(\exp(it)z).$ Choose a covering of
$S^1$ by intervals $R_{t_i}I,1\leq n \leq n$ and smooth functions
$f_i$ with support in $R_{t_i}I,1\leq n \leq n$ such that
$\sum_{1\leq i\leq n} f_i= \frac{1}{z},$ then
$$
0=\langle \xi, \sum_{1\leq i\leq n} V(\psi,f_i)\Omega\rangle=
\langle \xi, \psi\rangle, \forall \psi \in H_\lambda^0.
$$
Since $H_\lambda^0$ is dense in $H_\lambda,$ we conclude that
$\xi=0$ and the lemma is proved.
\endproof
%\qed
Recall $H_\B= \oplus_\la H_\a$ as representations of $LG_k$ or
$\A_{G_k}.$ The lowest energy space of $H_\a$, denoted by $H_\la(0)$
is a highest weight module of $G$ with weight $\la.$  The vertex
operator
$$
V(\psi,z): H_\la(0)\rightarrow \End(H_\B^0)[[z,z^{-1}]]
$$
is a primary vertex operator for $\hat g$ with highest weight
$\lambda$(cf. \cite{[KT]} and \cite{[FZ]}). By a slightly abuse of
notations we write such operator as $V(\la)=\sum_m V(\la)_m
z^{-m-1}.$

\begin{definition}\label{vv}
We define $V(\la)V(\mu)H_0$ to be the linear span of
$V(\la)_mV(\mu)_nH_0^0, \forall n,m.$
\end{definition}
Note that by definition $V(\la)V(\mu)H_0$ is a $\hat g$ submodule of
$H_\B,$ and if $V(\la)V(\mu)H_0\supset H_0^0,$ then
$\la=\bar\mu.$\par

The weight $1$ element in $H_\B$ is a Lie algebra isomorphic to
$h_{\Bbb C}$ and will be identified with $h_{\Bbb C}.$ It has a
subspace isomorphic to $g_{\Bbb C}.$ The vertex operator associated
with $h, V(h,z)$ is usually written as $h(z),$  similarly we write
$V(h,f)$ as $h(f).$ $h(f)$  are skew adjoint unbounded operators if
$f=f^*.$ The orthogonal complement of $g_{\Bbb C}$ in $h_{\Bbb C},$
denoted by $h_{\Bbb C}\ominus g_{\Bbb C}$ is a direct sum of
$H_\lambda(0)$ with $h_\la=1.$

\begin{lemma}\label{inc}
(1) Let $T_\la,T_\mu$ be as in Lemma \ref{charge}. Then

$$T_\lambda^*T_\mu^*H_0 \subset \overline{V(\lambda)V(\mu)H_0};
$$
(2) If $E(T_\nu T_\lambda^*T_\mu^*)\neq 0,$ then $H_\nu \subset
\overline{V(\lambda)V(\mu)H_0}; $
\end{lemma}
\proof Ad (1): We choose interval $I_1$ which is disjoint from $I.$
By Lemma \ref{RS} it is sufficient to check that for all smooth $f$
with support in $I_1$ and $\psi\in H_\mu^0,$
$$
T_\lambda^* V(\psi,f)\Omega \in \overline{V(\lambda)V(\mu)H_0}.
$$
By choosing $H$ trivial in Prop. 2.3 of \cite{X-coset}, we know that
$V(\psi,f)$ is affiliated with $\B(I_1),$ and by locality we have
$T_\lambda^* V(\psi,f)\Omega = V(\psi,f)T_\lambda^*\Omega.$ Since
$V(\mu)V(\lambda)H_0$ is a $\hat g$ module, the orthogonal
complement of $V(\mu)V(\lambda)H_0$ is a direct of irreducible $\hat
g$ module.

Now suppose that $\xi\in H^0_\B$ is orthogonal to
$V(\mu)V(\lambda)H_0.$ Choose $\chi_n\in H_\lambda^0$ such that
$\chi_n\rightarrow T_\lambda^*\Omega$ in norm. Then by Lemma 1 of
\cite{X-coset}
$$
\langle V(\psi,f)\chi_n,\xi\rangle = 0= \langle \chi_n, V(\psi,f)^*
\xi\rangle.
$$
Now let $n$ go to infinity we have
$$
\langle T_\lambda^*\Omega,V(\psi,f)^* \xi\rangle=0= \langle
V(\psi,f)T_\lambda^*\Omega, \xi\rangle,
$$
and (1) is proved. \par

Ad (2): Since $T_\la^*T_\mu^*= \sum_\nu T_\nu^* E(T_\nu
T_\la^*T_\mu^*),$ it follows that from (1) if $E(T_\nu
T_\la^*T_\mu^*)\neq 0,$ then $E(T_\nu T_\la^*T_\mu^*)^*\in \A(I)$ is
an isometry up to non-zero constant, and so
$H_\nu=T_\nu^*\overline{\A(I)\Omega} \subset
T_\la^*T_\mu^*\overline{\A(I)\Omega} \subset
\overline{V(\lambda)V(\mu)H_0}. $
\endproof
\begin{lemma}\label{wt1}
Suppose that $H_\lambda \in H_\C$ with $h_\lambda=1.$ Let $\psi\in
H_\lambda \cap (h\ominus g),$ and $f=f^*$ a smooth function with
support in $I.$ Then $\exp(V(\psi,f))\in \C.$
\end{lemma}
\proof Let $E_\C:\B\rightarrow \C$ be the conditional expectation
which is implemented by the projection $P_\C$ on $H_\B$ with range
$H_\C.$ We first show that
$$P_\C V(\psi,f)\exp(t V(\psi,f))\Omega = V(\psi,f)P_\C\exp(t
V(\psi,f))\Omega.
$$
For any $b\in \B(I')$ we have
\begin{align*}
\langle P_\C V(\psi,f)\exp(t V(\psi,f))\Omega,b\Omega \rangle
&=\langle V(\psi,f)\exp(t V(\psi,f))\Omega, E_\C(b)\Omega\rangle \\
&=-\langle\exp(t V(\psi,f))\Omega, V(\psi,f) E_\C(b)\Omega\rangle.
\end{align*}
Since by Prop. 2.3 of \cite{X-coset} $V(\psi,f)$ is skew self
adjoint and is  affiliated with $\B(I),$ and note that
$V(\psi,f)\Omega \in H_\lambda \subset H_\C,$ it follows that
\begin{align*}
\langle\exp(t V(\psi,f))\Omega, V(\psi,f) E_\C(b)\Omega\rangle &=
\langle\exp(t V(\psi,f))\Omega,  E_\C(b)V(\psi,f)\Omega\rangle\\
&=\langle\exp(t V(\psi,f))\Omega,  P_\C bV(\psi,f)\Omega\rangle\\
&=\langle P_\C \exp(t V(\psi,f))\Omega,  V(\psi,f)b\Omega\rangle.
\end{align*}

By (2) of Lemma 4 in \cite{X-coset} $P_\C \exp(t V(\psi,f))\Omega\in
H(\infty),$ and $H(\infty)$ is in the domain of skew self adjoint
operator $V(\psi,f).$ It follows that
$$
\langle P_\C \exp(t V(\psi,f))\Omega,  V(\psi,f)b\Omega\rangle=
-\langle V(\psi,f)P_\C \exp(t V(\psi,f))\Omega, b\Omega\rangle,
$$
and we have shown that
$$
\langle P_\C V(\psi,f)\exp(t V(\psi,f))\Omega,b\Omega
\rangle=\langle V(\psi,f)P_\C \exp(t V(\psi,f))\Omega,
b\Omega\rangle.
$$
By Reeh-Schleder Th. we have shown that $$P_\C V(\psi,f)\exp(t
V(\psi,f))\Omega = V(\psi,f)P_\C\exp(t V(\psi,f))\Omega.
$$

Set $F(t):= \langle P_\C\exp(t V(\psi,f))\Omega, \exp(t
V(\psi,f))\Omega\rangle.$ Then $F(0)=1$ and
\begin{align*}
F'(t)&= \langle V(\psi,f)\exp(t V(\psi,f))\Omega,P_\C\exp(t
V(\psi,f))\Omega\rangle +\\
&\langle \exp(t V(\psi,f))\Omega,P_\C V(\psi,f) \exp(t
V(\psi,f))\Omega\rangle =0
\end{align*}
where we have used $$P_\C V(\psi,f)\exp(t V(\psi,f))\Omega =
V(\psi,f)P_\C \exp(t V(\psi,f))\Omega
$$ and $P_\C \exp(t
V(\psi,f))\Omega\in H(\infty),$ and $H(\infty)$ is in the domain of
skew  self adjoint operator $V(\psi,f).$  It follows that $F(t)= 1$
and we conclude that
$$\exp(t V(\psi,f))\Omega\in H_\C$$ which proves our lemma.
\endproof
The following uses an analogue of VOA statement that weight $1$
space has a Lie algebra structure.
\begin{lemma}\label{chasweight1}
If $H_\lambda\subset H_\C$ with $h_\lambda=1,$ then $\C=\B.$
\end{lemma}
\proof By Lemma \ref{wt1} for any  $\psi\in H_\lambda \cap h,$ and
$f=f^*$ a smooth function with support in $I,$ we have
$\exp(V(\psi,f))\in \C.$  Since the conformal inclusions are
maximal, it follows that the Lie algebra generated by $g$ and $\psi$
is in fact Lie algebra $h.$ By Lie's formula, if $\exp(i
V(\psi_j,f_j))\in \C(I),j=1,2$ then
$$
((\exp(V(\psi_1,f_1)/n)\exp(V(\psi_2,f_2)/n)(\exp(-V(\psi_1,f_1)/n)\exp(-V(\psi_2,f_2)/n))^{n^2}$$
converges strongly to
$$
\exp([V(\psi_1,f_1),V(\psi_1,f_1)]).$$ On the other hand

$$[V(\psi_1,f_1),V(\psi_2,f_1)]= V([\psi_1,\psi_2],fg) + \langle
\psi_1,\psi_2\rangle \int_{S^1} f_1f_2 dz/z.
$$
It follows that for any $\psi\in h,$ and smooth functions $f=f^*$
with support in $I$ we have that
$$ \exp(V(\psi,f))\in \C(I).$$ Since $\B(I)$ is generated as a von
Neumann algebra by such elements, we have shown that $\C=\B$.
\endproof

\subsection{Induction of the adjoint representation}
The following is a key observation in this section, and is already
implicitly contained in (3) of Lemma 2.33 in \cite{X-adv}.
\begin{proposition}\label{key}
Suppose that $\C$ contains no weight $1$ element except those in
$\A,$ then $a_{v_0}^\C$ is irreducible.
\end{proposition}
\proof By Lemma \ref{fusionrule} we  have \[ [\ad^2]= [1]+ 2[\ad] +
[{(2,0,...,0,2)}] + [{(0,1,0,...,1,0)}]+ [{(0,1,0,...,0,2)}] +
[{(2,0,...,0,1,0)}]
\]
By computing the conformal dimensions of the descendants of $\ad^2$
using equation (\ref{cdim}) we have
\[
h_{(2,0,...,0,2)}= \frac{2+2n}{k+n},
h_{(0,1,...,0,2)}=h_{(2,0,...,1,0)} =\frac{2n}{k+n},
h_{(0,1,...,1,0)}= \frac{2n-2}{k+n} \]

Hence if $\C$ contains no weight $1$ element except those in $\A,$
then
$$
\langle a_{v_0}^\C, a_{v_0}^\C \rangle =\langle H_\C,v_0 v_0\rangle
=1$$ where recall that we use $H_\C$ to denote the restriction of
the vacuum representation of $\C$ to $\A,$ and the proposition is
proved.
\endproof
\begin{lemma}\label{localwithv0}
Suppose that $\epsilon(\lambda,v_0)\epsilon(v_0,\lambda)=1,$ then
$\lambda=\omega^i$ for some $0\leq i\leq n.$
\end{lemma}
\proof By definition we have
$$
\frac{S_{v_0\lambda}}{S_{1\lambda}}= d_{v_0}.
$$

From $[v][\bar v]=[1]+[v_0]$ we have

$$
\frac{S_{v\lambda}}{S_{1\lambda}}\frac{S_{\bar
v\lambda}}{S_{1\lambda}}=d_{v_0}+1 \leq d_vd_{\bar v} =d_{v_0}+1
$$
It follows that we must have
$|\frac{S_{v\lambda}}{S_{1\lambda}}|=d_v.$ For any positive integer
$k,$ suppose that $[v^k]=\sum_{\mu} m_\mu [\mu],$ then we have
$$
d_v^n= |\sum_{\mu} m_\mu \frac{S_{\mu\lambda}}{S_{1\lambda}}|\leq
\sum_{\mu} m_\mu |\frac{S_{\mu\lambda}}{S_{1\lambda}}|\leq
\sum_{\mu} m_\mu d_\mu = d_v^n.
$$
It follows that $|\frac{S_{\mu\lambda}}{S_{1\lambda}}|=d_\mu,
\forall \mu\prec v^k.$ Since every irrep of $\A$ occurs in some
$v^k,$ it follows that we must have
$|\frac{S_{\mu\lambda}}{S_{1\lambda}}|=d_\mu, \forall \mu.$ Square
both sides and sum over $\mu,$ we have proved that $d_\lambda=1,$
and hence the Lemma. \endproof

\subsection{List of intermediate subnets from conformal inclusions}
\begin{theorem}\label{c123}

(1) For the subnet $\A\subset \B$ corresponding to conformal
inclusions in \ref{c3}, when $n$ is odd (resp. even) the
intermediate subnet $\C$ are in one to one correspondence with the
abelian subgroup $\mathbb{Z}_n$ (resp. $\mathbb{Z}_{n/2}$) generated
by $\omega, $ (resp. $\omega^{2}$) i.e., if $\omega^i, ik=n$ (resp.
$\omega^{2i},2ik=n$) is a generator of this subgroup, then the
spectrum of $\C$ is  $H_\C=\sum_{1\leq j\leq k} H_{\omega^{ij}}$
(resp. $H_\C=\sum_{1\leq j\leq k} H_{\omega^{2ij}}$);\par

(2): For the subnet $\A\subset \B$ corresponding to conformal
inclusions in \ref{c1}, \ref{c2}, when $n$ is odd there is no
intermediate subnet. When $n=2m$ is even, the only nontrivial
intermediate subnet $\C$ is a $\mathbb{Z}_2$ extension of $\A$ by
the simple current $\omega^{m},$ i.e., the spectrum is $H_\C= H_0+
H_{\omega^m};$\par
\end{theorem}
\proof Ad (1): By Lemma \ref{chasweight1} we can assume that $\C$
has no weight $1$ elements besides those of $\A.$ By Prop. \ref{key}
we know that $a_{v_0}^\C$ is irreducible.  In the case of conformal
inclusions in \ref{c3}, since the vector representation of $LH,$
when restricting to $\A,$ contains the adjoint representation, it
follows from (4) of \ref{xua} that $a_{v_0}^\C$ must contain a DHR
representation of $\C.$ Since $a_{v_0}^\C$ is irreducible, it
follows that $a_{v_0}^\C$ is a DHR representation of $\C$, i.e.,
$[a_{v_0}^\C]=[\tilde a_{v_0}^\C].$ By Lemma \ref{loc2} we must have
for any $\lambda\in H_\C,$
$\epsilon(\lambda,v_0)\epsilon(v_0,\lambda)=1.$ By Lemma
\ref{localwithv0} we conclude that $\lambda=\omega^i$ for some
$0\leq i\leq n.$ (1) now follows easily by inspection of the
spectrum of $\A\subset \B$ in \cite{ABI}. \par
Ad (2): In the case
of conformal inclusions in \ref{c1} (resp. \ref{c2}) , we note that
the vector representation of $LH,$ when restricting to $\A,$
contains the antisymmetric representation $(0,1,0,...,0)$ (resp.
symmetric representation $(2,0,0,... 0)$ of $\A$).) Since
$$
\langle a_{v^2}^\C, a_{v^2}^\C\rangle = \langle
[a_{(2,0,...,0)}^\C]+[a_{(0,1,0,...,0)}^\C],
[a_{(0,2,0,...,0)}^\C]+[a_{(1,1,0,...,0)}^\C]\rangle =\langle
a_{v\bar v}^\C, a_{v\bar v}^\C\rangle =2,
$$
it follows that both $a_{(2,0,...,0)}^\C$ and $a_{(0,1,0,...,0)}^\C$
are irreducible. Hence as in the proof of (1), for  the case of
conformal inclusions in \ref{c1} (resp. \ref{c2}),
$a_{(0,1,0,...,0)}^\C$ (resp. $a_{(2,0,...,0)}^\C$) are DHR
representations of $\C.$ It follows that if $\lambda\in H_\C$, then
$\epsilon(\lambda,(0,1,0,...,0))\epsilon((0,1,0,...,0),\lambda)=1$
(resp.$\epsilon(\lambda,(2,0,...,0))\epsilon((2,0,...,0),\lambda)=1$).
Similarly
$\epsilon(\lambda,(0,0,...,0,1,0))\epsilon((0,0,...,1,0),\lambda)=1$
(resp.$\epsilon(\lambda,(0,0,...,2))\epsilon((0,0,...,2),\lambda)=1$).

Since by Lemma \ref{fusionv} $v_0$ appears in the product of
$(0,1,0,...,0)$ (resp. $(2,0,...,0)$)and its conjugate, it follows
that $\epsilon(\lambda,v_0)\epsilon(v_0,\lambda)=1.$ By Lemma
\ref{localwithv0} we conclude that $\lambda =\omega^i, 1\leq i\leq
n,$ and (2) follows by inspection of the spectrum of $\A\subset \B$
as given in \cite{LL}.
\endproof
We note that the same idea in the proof of Theorem above gives a
proof of the following:
\begin{corollary}
Suppose $\A_{SU(n)_k}\subset \C, n\neq n,n\pm 2$, and there is a
representation of $\C$, when restricting to $\A_{SU(n)_k},$ contains
$v_0.$ Then $\C$ is an extension by simple currents.
\end{corollary}

\begin{remark} Since conformal inclusion $SU(2)_{10}\subset Spin
(5)_1$ is not a simple current extension,  this example shows that
the condition in the above corollary is necessary. In fact in this
case the adjoint representation of $\A_{SU(2)_{10}}$ does not appear
in the restriction of any irreps of $\A_{Spin(5)_1}.$
\end{remark}

\endremark

\begin{theorem}\label{normalcase}
For the subnet $\A\subset \B$ associated with \ref{c4}, let
$(n,m)=p, n=n_1p, m=m_1p. $ Then  the intermediate subnets $\C$ are
in one to one correspondence with the subgroup of $\Bbb {Z}_p,$
i.e.,each such $\C$ has  spectrum $H_\C=\sum_{0\leq l\leq k_2}
H_{(\omega^{n_1k_1l},\dot{\omega}^{m_1k_1l}})$ with $k_1k_2=p,$
where we use $\dot{\la}$ to denote the highest weights of $SU(m)_n.$
\end{theorem}
\proof Since $\A_{SU(n)_m}\subset \B$ is normal (cf. \S4 of
\cite{X-m}), it follows that for each $(\lambda,\dot\lambda)\in
H_\C,$ we must have $ [a_\lambda^\C]= [a_{\dot\lambda}^\C], $ and
$\lambda\rightarrow a_\lambda^\C$ is a ring isomorphism. So if
$(\lambda_i,\dot\la_i)\in H_\C, i=1,2,$ and

$\la_3\prec \la_1\la_2,$ then $[a_{\la_3}^\C]\prec
[a^\C_{\la_1\la_2}]=[a^\C_{\dot\la_1\dot\la_2}],$ it follows there
must be a $\dot\la_3$ such that $[a_{\la_3}^\C]=[a_{\dot\la_3}^\C]$
and $(\lambda_3,\dot\lambda_3)\in H_\C.$

It follows that $\C$ are in one to one correspondence with the set
$R_\C$ of $\lambda$ with color zero mod $n$ which are closed under
conjugation and fusion product. If $R_\C$ contains any $\la$ with
$d_\la\neq 1,$ by Lemma \ref{fusionv} we have $v_0\in R,$ and it
follows that $R$ contains all  $\lambda$ with color zero mod $n,$ in
which case $\C=\B.$ Now assume that $d_\la=1$ if $\la\in R_\C.$ The
$R_\C$ must be a subgroup generated by $\omega^{n_1k_1},k_1k_2=p .$
Since the color of elements in $R_\C$ is zero mod $n,$  our Theorem
follows.
\endproof
By checking the list of intermediate subnets from  Th. \ref{c123}
and Th. \ref{normalcase} we immediately have:
\begin{corollary}
Conjecture \ref{netwall} is true for $\A_{G_K}\subset \A_{H_1}$
where $G_K\subset H_1$ are conformal inclusions in
\ref{c1},\ref{c2},\ref{c3} and \ref{c4}.
\end{corollary}

For the conformal inclusions $G_k\subset H_1,$ we write $V_\A$
(resp. $V_\B$) the VOA (cf. \cite{[FZ]}) associate with affine $\hat
g$ at level $k$ (resp. affine $\hat h$ at level $1$) We have natural
inclusion $V_\A\subset V_\B.$ We are interested in VOA $V_\C$ such
that $V_\A\subset V_\C\subset V_\B.$ We say a VOA is {\it simple} it
is irreducible as a representation over itself. Note that $V_\C$
will be direct sum of $\hat g$ modules, and we can write
$V_\C=\oplus H_\lambda$ and we refer to those $\la$ which appear in
$V_\C$ as the {\it spectrum} of $V_\C.$

\begin{proposition}\label{ctildec}
For any simple VOA $V_\C$ such that $V_\A\subset V_\C\subset V_\B,$
there corresponds a unique intermediate subnet $\C$ such the
spectrum of $\A\subset \C$ is the same as the spectrum of
$V_\A\subset V_\C.$
\end{proposition}
\proof Fix an interval $I.$ For each $\lambda $ in the spectrum of
$V_\A\subset V_\C,$ denote by $T_\lambda$ be as in Lemma
\ref{charge}. By Lemma \ref{inc}, it follows that if $E(T_\nu
T_\lambda^*T_\mu^*)\neq 0$ where $\lambda,\mu$ are in the spectrum
of $V_\A\subset V_\C,$ then $\nu$ is also in the spectrum of
$V_\A\subset V_\C.$  If $\lambda$ is in the spectrum of $V_\A\subset
V_\C$ but $\bar\lambda $ is not, then by the remark after Definition
\ref{vv} the action of $V_\C$ on $H_\lambda$ will span an invariant
subspace of $H_\C$ which does not contain $H_0,$ contradicting our
assumption that $\C$ is simple. By (2) of  Lemma \ref{charge}, $\A,
T_\lambda$ where $\lambda$ is in the spectrum of $V_\A\subset V_\C$
generate an intermediate subnet $\C$ with its spectrum the same as
the spectrum of $V_\A\subset V_\C.$
\endproof

The above proposition immediately implies the following theorem:
\begin{theorem}\label{voastatement}
The set of intermediate simple VOA $V_\C$  in $V_\A\subset V_\B$ for
conformal inclusions \ref{c1},\ref{c2},\ref{c3} and \ref{c4} are in
one to one correspondence with the set of intermediate subnets $\C$
of $\A\subset \B$ with the same spectrum as given in Th. \ref{c123}
and Th. \ref{normalcase}.
\end{theorem}
\begin{remark}
We note that the simple intermediate VOAs in the above theorem are
simple current extensions of affine VOAs, and they are well
understood in VOA literature (cf.\cite{DL}).
\end{remark}

\section{Verifying Conjecture \ref{wall} for Jones-Wassermann subfactors}

In this section we extend the results in Cor. 5.23 of \cite{X-adv}.

Let $\la$ be an irreducible representation of $\A_{SU(n)_k}$
localized on $I, M:=\A_{SU(n)_k}(I). $ Suppose $\la=c_1c_2$ where
$c_i\in \End(M),i=1,2,c_1(M)$ is an intermediate subfactor of
$\la(M)\subset M.$ We note that $c_1\bar c_1\prec \la\bar\la.$ We
say the intermediate subfactor $c_1(M)$ is of {\it abelian type} if
$[c_1\bar c_1]=\sum_{1\leq i\leq j} [\omega^{ij_1}], jj_1=n.$ The
following Lemma appears as Lemma 5.22 in the correction of proof of
\cite{X-adv} and we include its proof: \blem\label{c12} Assume that
$Z_{1\mu}^{c_1}=\delta_{1\mu}, \forall \mu$ where $Z^{c_1}$ is
defined as in Definition \ref{zmatrix}. Then $\langle
c_1c_2,c_1c_2\rangle = \langle c_1\bar c_1,\bar c_2
c_2\rangle.$\elem

\prf By \S2 of \cite{Ganv} we have
$Z_{\mu_1\mu_2}^{c_1}=\delta_{\mu_1\tau(\mu_2)}$ where
$\mu\rightarrow \tau(\mu)$ is an order two automorphism of fusion
algebra. It follows that $[\tilde a_\mu]=[a_\tau(\mu)],$ and by
\cite{BEK2}  irreducible sectors of $\bar c_1 \nu c_1$ are of the
form $a_\mu, \forall \mu.$ Since
$$
\langle c_2\bar c_2, a_\mu\rangle= \langle c_2, a_\mu c_2\rangle=
\langle c_2, c_2\mu\rangle=\langle \bar c_2 c_2,
a_\mu\rangle=\langle a_{\bar c_2 c_2}, a_\mu\rangle,
$$
we conclude that $[c_2\bar c_2]=[a_{\bar c_2 c_2}],$ and
$$
\langle c_1\bar c_1, \bar c_2 c_2\rangle= \langle c_1, \bar c_2 c_2
c_1 \rangle= \langle c_1, c_1 a_{\bar c_2 c_2}\rangle=\langle c_1,
c_1 c_2 \bar c_2\rangle=\langle { c_1 c_2}, c_1 c_2\rangle
$$
\qed

\begin{theorem}\label{max} (1) Suppose that $k\neq n-2,n+2,n.$  then $\la$ is
maximal iff there is no $ 1\leq i\leq n-1$ such that
$[\omega^i\la]=[\la];$ \par (2) When $\la$ is not maximal, the
maximal intermediate subfactor is either abelian type or at most one
given by $c_1(M)$ with $\la=c_1c_2, [\bar
c_2][c_2]=[1]+[\omega^m],n=2m.$

\end{theorem}

\prf Ad (1): (1) is Cor. 5.23 in \cite{X-adv}. We include its proof
which will be modified in our proof of (2).\par

When $k=1$  the Cor. is obvious. By Lemma 2.33 of \cite{X-adv} we
can assume that $k\geq 2$ and $d_\ad>1.$
 As in the proof of Cor. 5.21 in \cite{X-adv}, $\la$ is maximal implies that
there is no $ 1\leq i\leq n-1$ such that $[\omega^i\la]=[\la].$ Now
suppose that there is no $ 1\leq i\leq n-1$ such that
$[\omega^i\la]=[\la].$ If $S_{v\la}\neq 0$, then $\la$ is maximal by
Cor. 5.20 of \cite{X-adv}. If $k=2,$ the $S$ matrix elements are
equal to that of $S$ matrix elements for $SU(2)_n$ up to phase
factors, and it follows easily that $S_{v\la}\neq 0$ if there is no
$ 1\leq i\leq n-1$ such that $[\omega^i\la]=[\la].$

Suppose that $k \geq 3, S_{v\la}= 0.$ Since $[v\bar v]=[1]+[\ad]$ we
have $S_{\ad\la}= -S_{1\la}\neq 0.$ Assume that $M_1$ is an
intermediate subfactor between $\la(M)$ and $M$, and $\la=c_1c_2$
with $c_1(M)=M_1$ and $c_1=c_1'c_1''$ as in Prop. \ref{loclr}. Apply
Lemma 2.20 of \cite{X-adv}  we have $\langle a^{c_1'}_{\ad}, \tilde
 a^{c_1'}_{\ad}\ra \geq 1.$ By Lemma 2.33 of \cite{X-adv} we must have
 $[ a^{c_1'}_{\ad}]=[ \tilde a^{c_1'}_{\ad}]$ and by Lemma
 2.36 of \cite{X-adv}  $[c_1'\bar c_1']= \sum_{1\leq j\leq n/j_1} [\om^ {j j_1}].$
 By Frobenius reciprocity we have $[\om^{j_1}c_1']=[c_1'.]$
Since $\la=c_1'c_1''c_2,$ $[\om^{j_1}\la]=[\la]$, and by assumption
$j_1=n$ and $[c_1'\bar c_1']=[1]$. By  Prop. \ref{loclr} we must
have   $Z_{\mu1}^{c_1}= \delta_{\mu1}, \forall \mu.$ By \S2 of
\cite{Ganv} we have $Z_{\mu_1\mu_2}^{c_1}=\delta_{\mu_1\tau(\mu_2)}$
where $\tau(\mu)=\omega^{m \col(\mu)} \mu  $ or $\tau(\mu)=\omega^{m
\col(\mu)} \bar\mu,  m\geq 0.$ We claim that in fact
$[\omega^m]=[1]$ and $\tau(\mu)=\mu.$ First we show that
$\tau(\mu)=\omega^{m \col(\mu)} \mu.  $  If instead
$\tau(\mu)=\omega^{m \col(\mu)} \bar\mu,$
 since $k\geq 3,$
$\tau((0,1,0,...,0))\neq (0,1,0,...,0),$ by Lemma 2.20 of
\cite{X-adv} we must have $S_{\la(0,1,0,...,0)}=0.$ From the fusion
rule
$$
[(0,1,0,...,0)(0,0,...,0,2)]=[(0,1,0,...,0,2)]+[v_0]
$$
we must have $S_{\la(0,1,0,...,0,2)}\neq 0.$ By Lemma 2.20 of
\cite{X-adv}  we must have
$\tau((0,1,0,...,0,2))=(0,1,0,...,0,2)=(2,0,0,...,1,0),$ a
contradiction.  So we conclude that $\tau(\mu)=\omega^{m \col(\mu)}
\mu,\forall \mu$. It follows that $[\tilde a_\mu]=[a_{\omega^{m
\col(\mu)}}a_\mu],$ and in particular $[\tilde
a_v]=[a_{\omega^{m}}a_v]$. So we have
$$
[\omega^m v c_1]=[c_1 \tilde a_v]=[c_1 a_v]=[vc_1],
$$
and similarly $[c_2 \omega^{-m}\bar v ]=[c_2\bar v]$. If
$[\omega^m]\neq [1],$ by our assumption on $\la$ we have
$\omega^m\not\prec c_1\bar c_1, \omega^m\not\prec \bar c_2 c_2 .$ On
the other hand we have
$$
\langle \bar v \omega^m v,c_1\bar c_1\rangle \geq 1, \langle \bar v
\omega^m v,\bar c_2 c_2\rangle \geq 1
$$
It follows that $\omega^m v_0\prec c_1\bar c_1, \omega^m v_0\prec
\bar c_2 c_2,$ and $\langle c_1\bar c_1, \bar c_2 c_2\rangle \geq
2.$ By Lemma \ref{c12} we conclude that $\la=c_1c_2$ is not
irreducible, contradicting our assumption. Hence $[\omega^m]=[1]$
and $Z_{\mu_1\mu_2}=\delta_{\mu_1\mu_2}.$  The rest of the proof now
follows in exactly the same way as in the proof of Prop. 5.20 of
\cite{X-adv}.
\par
Ad (2): As in the proof of (1) we assume that $\la=c_1c_2$ with
$c_1(M)$ a nontrivial maximal intermediate subfactor, and
$c_1=c_1'c_1''.$ By our assumption we must have $c_1=c_1'$ if
$[c_1']\neq [1].$ In this case as in the proof of (1) above we must
have $[c_1\bar c_1]= \sum_{1\leq j\leq n/j_1} [\om^ {j j_1}].$\par
Now suppose that $[c_1']=[1].$ Then as above we have
$Z_{\mu\gamma}^{c_1}=\delta_{\mu,\omega^{\mu\col(\gamma)}\gamma}.$
By Corollary 3.14 of \cite{X-adv} we can find $c\prec \mu c_1$ for
some $\mu$  such that
$$[c\bar c]=\sum_{1\leq i\leq p} \omega^{li}, pl=n.$$
Since $Z^c=Z^{c_1},$ it follows the left local support of $c$ is
trivial.

If $[\omega^m]=[1],$ then we are as in the end of proof of (1), and
in that case $[c_1]=[1],$ contradicting our assumption that $c_1(M)$
a nontrivial maximal intermediate subfactor. So $[\omega^m]\neq[1].$
If $\omega^m\prec c_1\bar c_1,$ then by maximality of $c_(M)$ we
have $[c_1\bar c_1]=\sum_{i}[\omega^{qi}],$  i.e., $c_1(M)$ comes
from abelian part of $\la.$ Now assume that $\omega^m \not \prec
c_1\bar c_1,$ then as in (1) we must have $\omega^m\prec \bar
c_2c_2.$

Since $c\prec \mu c_1$ for some $\mu,$ we have $c\bar c\prec \mu
c_1\bar c_1\bar\mu\prec \mu\lambda\bar\lambda\bar\mu,$ so
$\col(\omega^{li})=0 \ \mod \ n,$ so we have $n|li.$

Similarly since $\omega^m \prec \bar c_2 c_2,$ $\col(\omega^m)=0 \
\mod \ n.$ On the other hand since the map $\mu\rightarrow
\omega^{m\col(\mu)}\mu$ has order two, it follows that
$\omega^{2m}=1.$ So we must have $n=2m.$

From
$$h_{\omega^{li}}=\frac{kli}{n}\frac{n-li}{2},$$ it follows that
$h_{\omega^{li}}\in \Z$ if $i$ is even. These $\omega^{li}$ with $i$
even will generate local simple currents (cf Definition 2.3 and
Prop. 2.15 of \cite{X-am}), and it follows that the left local
support of $c$ is nontrivial if $li\neq 0\mod n$ for some even $i.$
So we conclude that $2l=0\mod n,$ and $[c\bar c]=[1]+[\omega^m].$
\par

Note that there are $\lambda_1,\lambda_2$ such that $c_1\prec
\lambda_1 c, \bar c_2\in \lambda_2 c.$ From $\langle c_1,\lambda_1
c\rangle=\langle c_1\bar c,\lambda_1\rangle \geq 1,$ we have
$d_{c_1}\sqrt{2} \geq d_{\lambda_1},$ and similarly $d_{\bar
c_2}\sqrt{2} \geq d_{\lambda_2}.$

Since $[\omega^m c_1]\neq [c_1],$ we have $[\lambda_1 c]\succ
[c_1]+[\omega^m c_1],$ and by computing statistical dimension
$[\lambda_1 c]=[c_1]+[\omega^m c_1],$ and $[c_1\bar c]=[\lambda_1].$

Similarly if $\lambda_2 c$ is not irreducible, we must have $[\bar
c_2 \bar c]=[\lambda_2].$ But we have
$$
\langle \bar c_2 \bar c,\bar c_2 \bar c\rangle = \langle \bar c_2
\bar c c c_2,1\rangle =\langle \bar c_2 a_{c\bar c} c_2,1 \rangle
=\langle  c\bar c \bar c_2c_2,1 \rangle\geq 2,
$$
where we have used $[\bar c c] = [a_{c\bar c}]$ and $\bar
c_2c_2\succ [1]+[\omega^m].$ From this we conclude that $[\bar
c_2]=[\lambda_2 c].$\par

We have $[\lambda]=[c_1c_2]=[c_1\bar c \bar
\lambda_2]=[\lambda_1\bar \lambda_2].$ By Lemma \ref{fusionv} , we
must have $\bar \lambda_2=\omega^i,[\lambda] =[c_1 \bar c
\omega^i].$

Now we check that the intermediate subfactor $c_1(M)$ is uniquely
fixed. Suppose there is an intermediate subfactor $f_1(M)$ such that
$[\lambda]=[f_1f_2]$ and $[\bar f_2 f_2]=[1]+[\omega^m].$ Then $\bar
c\bar f_2$ has statistical dimension two and decompose into sum of
two irreducible endomorphisms, it follows that there is an
automorphism $\alpha$ such that $[f_2]=[\alpha \bar c\omega^i],$ and
$[f_1\alpha ]=[c_1\beta]$ for some automorphism  $\beta\prec \bar c
c.$ By Cor. 2.4 of of \cite{X-ghj} the intermediate subfactor
$f_1(M)$ is determined by equivalence class $[f_1,f_2]$ with
equivalence relation $[f_1,f_2]\sim [f_1 \rho,\rho^{-1}f_2]$ where
$\rho$ is any automorphism. We have
$$[f_1,f_2]=[c_1\beta
\alpha^{-1},\alpha\bar c\omega^i]\sim [c_1\beta, \bar c\omega^i]\sim
[c_1,\beta^{-1}\bar c\omega^i]=[c_1,\bar c\omega^i].
$$

\endproof

\begin{corollary}\label{checkjw}
Suppose that $k\neq n-2,n+2,n.$ Then each irreducible representation
$\lambda$ of $\A_{SU(n)_k}$ verifies both maximal and minimal
version of Conjecture \ref{wall}.
\end{corollary}
\proof Since the dual of $\lambda$ is $\bar\lambda,$ it is
sufficient to verify the maximal version of Conjecture \ref{wall}.
We may assume that $d_\lambda > 1.$ By Lemma \ref{fusionv}
$$
[\lambda][\bar \lambda]=\sum_{1\leq i\leq p}[\omega^{iq}]+[v_0]+...
$$
where $pq=n,$ and $...$ are possible additional irreps. We note that
the set of maximal intermediate subfactors for $\lambda$ coming the
abelian part are bounded by $p-1,$ and by Th. \ref{max}, there is at
most one more maximal intermediate subfactor, and our corollary
follows.
\endproof
\begin{remark}
It will be interesting to remove the condition $k\neq n,n\pm2 $ in
Th.\ref{max} and Cor. \ref{checkjw}. This condition is used in the
proof of Th. \ref{max} to ensure that $a_{v_0}^{c_1}$ is
irreducible. One can remove this condition if one can find a
different way of proving that $a_{v_0}^{c_1}$ is irreducible.
\end{remark}

\end{document}